\documentclass[12pt]{article}

\usepackage[english]{babel}
\usepackage{amsmath,amsfonts,amssymb,latexsym,graphics,epsfig,url}
\usepackage{color}
\usepackage{amsthm}
\usepackage{graphicx}
\usepackage{soul}
\usepackage[utf8]{inputenc}
\usepackage{graphicx}
\usepackage{latexsym}
\usepackage{mathtools}
\usepackage{anyfontsize}

\newtheorem{theorem}{Theorem}[section]
\newtheorem{proposition}[theorem]{Proposition}

\newtheorem{example}[theorem]{Example}

\newtheorem{remark}[theorem]{Remark}

\DeclareMathOperator{\diag}{diag}

\DeclareMathOperator{\Irep}{Irep}
\DeclareMathOperator{\rank}{rank}
\DeclareMathOperator{\tr}{tr}
\DeclareMathOperator{\spec}{sp}

\def\c{\mbox{\boldmath $c$}}

\def\f{\mbox{\boldmath $f$}}

\def\0{\mbox{\boldmath $0$}}

\def\vec0{\mbox{\boldmath $0$}}

\def\A{\mbox{\boldmath $A$}}
\def\B{\mbox{\boldmath $B$}}

\def\G{\Gamma}

\def\O{\mbox{\boldmath $O$}}
\def\S{\mbox{\boldmath $S$}}

\def\T{\mbox{\boldmath $T$}}
\def\U{\mbox{\boldmath $U$}}

\DeclareMathOperator{\Z}{\mathbb{Z}}

\def\1{\mbox{\boldmath $1$}}

\def\c{\mbox{\boldmath $c$}}

\def\CC{\mathbb C}

\def\G{\Gamma}

\def\Z{\mathbb Z}

\begin{document}
\title{Spectra and eigenspaces of branched covers \\ of graphs induced by local factorisation
}

	\author{C. Dalf\'o$^a$, M. A. Fiol$^b$,  S. Pavl\'ikov\'a$^c$, and J. \v{S}ir\'a\v{n}$^d$\\
		\\
		{\small $^a$Departament. de Matem\`atica}\\ {\small  Universitat de Lleida, Igualada (Barcelona), Catalonia}\\
		{\small {\tt cristina.dalfo@udl.cat}}\\
		{\small $^{b}$Departament de Matem\`atiques}\\
  {\small Universitat Polit\`ecnica de Catalunya, Barcelona, Catalonia} \\
		{\small Barcelona Graduate School of Mathematics} \\
		{\small  Institut de Matem\`atiques de la UPC-BarcelonaTech (IMTech)}\\
		{\small {\tt miguel.angel.fiol@upc.edu}}\\
        {\small $^{c}$ Inst. of Information Engineering, Automation, and Math.} \\
   {\small Slovak University of Technology, Bratislava, Slovak Republic}\\
		{\small {\tt sona.pavlikova@stuba.sk} }\\
   {\small $^{d}$ Department of Mathematics and Descriptive Geometry}\\
   {\small Slovak University of Technology, Bratislava, Slovak Republic}\\
		{\small {\tt  jozef.siran@stuba.sk}}
 }

\date{}
\maketitle

\begin{abstract}
{\small
We consider lifting eigenvalues and eigenvectors of graphs to their {\em factored lifts}, which are branched covers formed from unbranched ones by what one can describe as a `local factorisation' induced by subgroups of vertex stabilisers. Such a description applies to all situations in which a group of automorphisms of a branched cover of a graph acts freely on the arc set of the cover. Branched covers with this property can be described by means of {\em combined voltage assignments} in a group, which generalise the (ordinary) 
voltage assignments known from the theory of regular covering
and are a particular case of an algebraic description of graph covers in Poto\v{c}nik and Toledo (2021).

With the help of group representations and certain matrices over complex group rings associated with the graphs to be branch-covered, we develop a method for determining the complete spectra of factored lifts, and derive a sufficient condition for lifting eigenvectors.}

\end{abstract}

\section{Introduction}

A prolific construction of new graphs from old relies on regular covers, stemming from the theory of covering spaces in algebraic topology and made popular among graph theorists through the monograph by Gross and Tucker \cite{gt87}. Informally, one starts with a `base graph' equipped with a `voltage assignment' on its arcs, with values in some group. The resulting `voltage graph' gives rise to a new graph that may be considered as a `lift' of the base graph. The vertex- and edge-set of the lift is a product of the vertex- and edge-set of the base graph with the voltage group, and incidence in the lift is defined in such a way that the voltage group acts freely on the vertex-set of the lift as a group of automorphisms. Such a lift is a regular covering space of the base graph, where the covering is induced via a graph homomorphism given simply by erasing the group coordinate. Conversely, if one has a graph with a given group of automorphisms acting freely on its vertices, then the graph arises as a regular lift as indicated, with the base graph isomorphic to the quotient of the original graph by the given group of automorphisms.

Prominent examples of graph covers include Cayley graphs, usually introduced as graphs admitting a group of automorphisms acting regularly on their vertex set, which is equivalent to arising as regular lifts of one-vertex graphs with loops and semi-edges. But constructions based on graph covers extend well beyond vertex-transitivity and range from graphs of given degree and diameter (or girth) and expanders to milestones such as the Map Color Theorem. In a nutshell, the advantages of covering constructions lie in the fact that, in a number of important situations, properties of a lift can be conveniently expressed in terms of properties of a base graph and a corresponding voltage assignment.

Despite numerous applications of coverings, relatively little attention has been paid in the past to the related `lifting'  of eigenvalues and eigenspaces of a base graph to its covers. A notable advancement in the determination of spectra and eigenspaces was made in  Dalf\'o, Fiol, Miller, Ryan, and \v{S}ir\'a\v{n} \cite{dfmrs19}, with adaptation to digraphs in Dalf\'o, Fiol,  and \v{S}ir\'a\v{n} \cite{dfs19}, generalisation to non-regular lifts of graphs in Dalf\'o, Fiol, Pavl\'ikov\'a, and \v{S}ir\'a\v{n} \cite{dfps21}, and further extension to the universal adjacency matrix of such lifts in \cite{dfps23} by the same authors.

Our aim is to further extend these advances by introducing the concept of a `factored lift,' arising by considering a free action of a subgroup of automorphisms on the {\em arc set} instead of the vertex set of a graph. Since in such a situation the free action of a subgroup on arcs is not necessarily free on vertices, it induces, in general, covers with `branch points' at vertices, and hence corresponds to branched coverings in algebraic topology. Such a viewpoint is also a critical special case of the recently developed general approach to coverings in Potocnik and Toledo \cite{pt21}, allowing an action of a completely arbitrary subgroup of automorphisms. Our factored lifts, independently introduced in terms of `overlifts' in 
Reyes, Dalf´o, Fiol, and Messegu\'e \cite{rdfm23}, also represent a significant generalisation of permutation voltage lifts (describing non-regular unbranched graph covers), with implications for broader theoretical and practical applications.

The fact that factored lifts may induce graph covers with branch points at vertices makes the study of lifts of spectra and eigenspaces way more complex, compared to covers without branch points. The latter can equivalently be characterised in terms of the so-called permutation voltage assignments, suitable also for the description of non-regular covers (those not necessarily induced by group actions). A complete characterisation of lifts of spectra and eigenspaces of graph covers given by permutation voltage assignments was given in Dalf\'o, Fiol, Pavl\'ikov\'a, and \v{S}ir\'a\v{n} \cite{dfps21}. Such a characterisation, however, appears to be elusive in the more general setting of factored lifts, where the troublemakers are lifts of `singular eigenspaces', those corresponding to zero eigenvalues. Building a generalisation of a theory of lifts of spectra and eigenspaces for branched covers of graphs and identifying obstacles for a complete characterisation (accompanied by illustrative examples) has been the main motivation for our work.

The structure of the paper is as follows. In the next section, we give definitions and a formal statement of the equivalence between factored lifts and quotients by a free action of an automorphism group of a graph on its arcs. Section \ref{sec:walks} is devoted to lifts of walks and their enumeration. Our main results on lifts of spectra and eigenvectors from base graphs to factored lifts are in Section \ref{sec:factored-lifts}. In Section \ref{sec:examples},  we illustrate our results on two examples, followed by concluding remarks in Section \ref{sec:rem}.


\section{Combined voltage assignments and factored lifts}\label{sec:comb-fact}

All graphs in this article are finite and undirected, but may contain loops, multiple edges, and semi-edges. Recall that a semi-edge is incident to just one vertex, with its other end `dangling' and not incident to any vertex. (If graphs are identified with one-dimensional CW-complexes, a semi-edge is a quotient of a complete graph $K_2$ by an automorphism of order $2$ interchanging its vertices.) For the sake of convenience, however, we identify every edge $e$ joining two vertices $u$ and $v$ (possibly, with $u=v$) with a pair of oppositely directed {\em arcs}, arising from $e$ by assigning it a direction (from $u$ to $v$, or from $v$ to $u$). If $a$ is an arc arising from such an edge $e$, its {\em reverse} arc is denoted $a^-$. A semi-edge gives rise to just one arc $a$, identified with the semi-edge itself, and in this exceptional case, we let $a=a^-$.

Let $\Gamma$ be a finite graph with vertex set $V=V(\Gamma)$ and arc set $A=A(\Gamma)$, and let $G$ be a finite group.  Let $\alpha:\ A\to G$ be a voltage assignment on $\Gamma$ in the usual sense (see Gross and Tucker \cite{gt87}), that is, with $\alpha$ satisfying $\alpha(a^-)=\alpha(a)^{-1}$ for any arc $a\in A$ and its reverse $a^-$. (It follows that the value of $\alpha$ on a semi-edge is an element of $G$ of order at most two.) Let $\omega$ be a function assigning to every vertex $u\in V$ a subgroup $G_u$ of $G$. The pair $(\alpha,\omega)$ will be called a {\em combined voltage assignment} on the {\em base graph} $\Gamma$, and the graph thus becomes a {\em combined voltage graph}, with {\em voltage group} $G$.

The {\em factored lift} $\Gamma^{(\alpha,\omega)}$ is a graph with vertex set $V^{(\alpha,\omega)}$ consisting of all pairs $(u,H)$, where $u\in V$ and $H\in G/G_u = \{hG_u\ |\ h\in G\}$, the set of {\em left} cosets of $G_u$ in $G$. The arc set $A^{(\alpha, \omega)}$ of the factored lift is defined as follows. For a pair $u,v\in V$ of adjacent vertices in the base graph $\Gamma$, let $\overrightarrow{uv}$ denote the set of arcs of $\Gamma$ emanating from $u$ and terminating at $v$, which may coincide with $u$; if $v$=$u$ the set $\overrightarrow{uv}$ is also assumed to contain all semi-edges at $u$. Let $a\in \overrightarrow{uv}$ be an arc carrying a voltage $\alpha(a)\in G$. Then, for every $h\in G$, the arc $a$ together with the element $h$ determine a unique arc $(a,h)\in A^{(\alpha, \omega)}$ in the factored lift, emanating from the vertex $(u,hG_u)\in V^{(\alpha,\omega)}$ and terminating at the vertex $(v,h\alpha(a)G_v)\in V^{(\alpha, \omega)}$. Equivalently, if $a\in \overrightarrow{uv}$ is an arc in the base graph $\Gamma$ carrying voltage $\alpha(a)\in G$ and if $H\in G/G_u$ and $K\in G/G_v$ are left cosets, then
\begin{equation}\label{eq:def}
(u,H)\xlongrightarrow[]{(a,h)} (v,K)\  {\rm for\ every} \  h\in H \ {\rm and\ each}\ K \ {\rm such\ that}\  h\alpha(a)\in K\, 
\end{equation}
where the symbol on the left-hand side of \eqref{eq:def} denotes the arc $(a,h)$ from $(u,H)$ to $(v,K)$ in the factored lift. 

It is important to note that the factored lift $\Gamma^{(\alpha,\omega)}$ can, and will, be considered to be an {\em undirected graph}. Indeed, if $a\in \overrightarrow{uv}$ is an arc of $\Gamma$ that is not a semi-edge, then, for every $h\in G$ the reverse arc $(a,h)^-$ to the arc $(a,h)$ emanating from $(u,hG_u)$ and terminating at $(v,h\alpha(a)G_v)$, may be defined by letting $(a,h)^- = (a^-,h\alpha(a))$. To see this, observe that since $a^-\in \overrightarrow{vu}$, in the factored lift the arc $(a^-,h\alpha(a))$ emanates from the vertex $(v,h\alpha(a)G_v)$ and terminates at the vertex $(u,h'G_u)$, where, by definition of arcs in $\Gamma^{(\alpha,\omega)}$, one has $h' = (h\alpha(a))\alpha(a^-) = h$, so that $(a,h)^-$ terminates at $(u,hG_u)$ as expected. For a semi-edge $a$ in $\Gamma$ at a vertex $u$, one has two possibilities: either $\alpha(a)=1$, in which case $(a,h)$ is a semi-edge at the vertex $(u,hG_u)$ in the factored lift for every $h\in G$, or else $\alpha(a)$ is an involution. In the latter case, the arguments in the first part of this paragraph apply, with $(a,h)$ being an edge from $(u,hG_u)$ to a different vertex $(u, h\alpha(a)G_u)$ if $\alpha(a)\not\in G_u$, and a loop at the vertex $(u,hG_u)$ if $\alpha(a)\in G_u$ but $\alpha(a)$ is an involution; in both cases, one has $(a,h)^- = (a^-,h\alpha(a))$, as before.

\begin{figure}[t]
	\centering
\includegraphics[width=14cm]{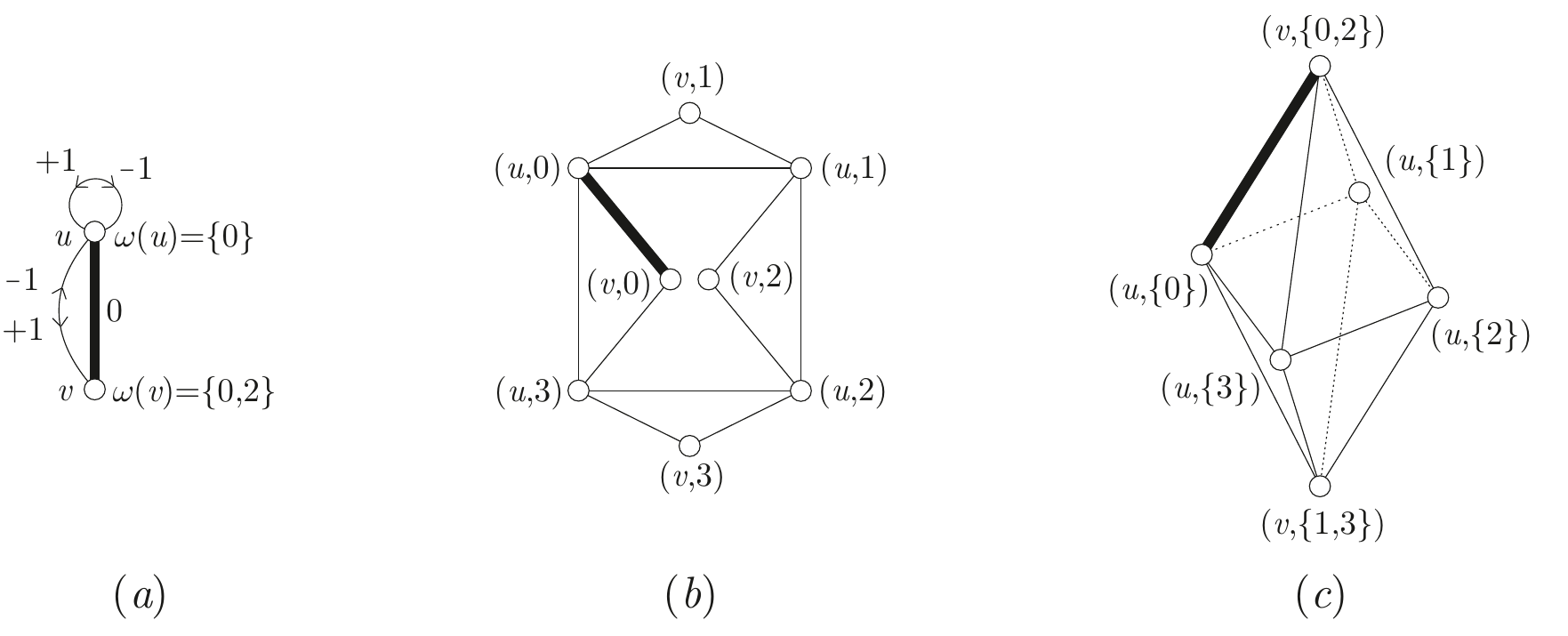}
	\caption{$(a)$ A combined voltage graph $\G$ on $\Z_4$. $(b)$ The ordinary lift of $\G$. $(c)$ The factored lift of $\G$.}
  \label{fig1}
\end{figure}

As a simple example, Figure \ref{fig1} shows a combined voltage graph on the cyclic group $\Z_4$ (on the left) and the resulting factored lift (the graph of an octahedron) on the right.

If the subgroup $G_u < G$ is trivial for every vertex $u$ of $\Gamma$, a factored lift reduces to an {\em ordinary} lift $\Gamma^{\alpha}$ with vertex set $V^{\alpha} = \{(u,h)\ |\ u\in V, h\in G\}$, in which for every arc $a \in \overrightarrow{uv}$ of $\Gamma$ and every $h\in G$ the pair $(a,h)$ is an arc in the ordinary lift $\Gamma^{\alpha}$ from the vertex $(u,h)$ to the vertex $(v,h\alpha(a))$, see Gross and Tucker \cite{gt87}. The way a factored lift $\Gamma^{(\alpha,\omega)}$ arises from an ordinary lift $\Gamma^{\alpha}$ should now be obvious. Indeed, for a general assignment $\omega:\ u\mapsto G_u$ and for every $u\in V$ one identifies {\em right} $G_u$-orbits in the fibre $\{(u,h)\ |\ h\in G\}$ of $\Gamma^{\alpha}$ to form vertices $(u,hG_u)$, but making no identification among the existing arcs. This process can be regarded as a `factorisation' induced by `local right actions' of the subgroups $G_u<G$ for $u\in V$; hence the term {\em factored lift}.

To illustrate this by the example of Figure \ref{fig1}, if all the vertices of the voltage graph are assigned the trivial group, the obtained (standard) lift graph is shown in Figure \ref{fig1} (centre). The factored lift on the right of Figure \ref{fig1} is then obtained by identifying the pairs $\{(v,0),(v,2)\}$ and $\{(v,1),(v,3)\}$ to a single vertex each, resulting in the graph of an octahedron.


The factorisation induces a graph epimorphism $f: \Gamma^{\alpha} \to \Gamma^{(\alpha,\omega)}$, given simply by $(u,h) \mapsto (u,H)$ for $H\in G/G_u$ such that $h\in H$ (we have deliberately replaced $g$ by $h$ in the notation for future use). A more detailed description of the action of $f$ is on the following diagram:

\begin{equation}\label{eq:diag}
\begin{array}{ccc}
(u,h) & {\stackrel{(a,h)}{\xlongrightarrow{\hspace*{2cm}} }} & (v,h\alpha(a)) \\
& & \\ \Big\downarrow {\scriptstyle {f}} & & \Big\downarrow {\scriptstyle {f}} \\
(u,hG_u) & {\stackrel{(a,h)}{\xlongrightarrow{\hspace*{2cm}} }} & (v,h\alpha(a)G_v)
\end{array}
\end{equation}

This factorisation results in $[G:G_u]$ vertices of the form $(u,H)$ for left cosets $H$ of $G_u$ in $G$, each as an $f$-image of the $|G_u|$ vertices $(u,h)$ for $h\in H$. Every vertex $(u,H)$ then has valency $|G_u|d_u$, where $d_u$ is the valency of $u$ in $\Gamma$. Moreover, the group $G$ acts on the factored lift $\Gamma^{(\alpha,\omega)}$ by left multiplication as a subgroup of automorphisms. The action is {\em free} on the arc set $A^{(\alpha,\omega)}$, and is transitive on every fibre $F(u)=\{(u,H)\in V^{(\alpha, \omega)}\ |\ H\in G/G_u\}$, with $G_u$ being the stabiliser of the vertex $(u,G_u)\in V^{(\alpha,\omega)}$ for every $u\in V$. In particular, this action of $G$ produces the base graph $\Gamma$ as a $G$-quotient of its factored lift; formally, it gives rise to a graph isomorphism $\Gamma \cong \Gamma^{(\alpha, \omega)}/G$.

We remark that combining {\em left} cosets of $G$ with {\em right} multiplication of elements of $G$ by $\alpha(a)$ in \eqref{eq:def} is essential for the algebra in the factorisation and morphisms to work.

Factored lifts turn out to be a special case of the {\em generalised voltage assignments and lifts}, introduced by Poto\v{c}nik and Toledo \cite{pt21}, where the objects assigned to arcs are also allowed to be cosets of the voltage group rather than individual elements of this group. In this generalisation, both vertices and arcs of a lift are `indexed' by subgroups of the voltage groups, subject to several technical conditions. The difference between our treatment and the (more general) one in \cite{pt21} is the use of right cosets and left multiplication by voltages in \cite{pt21} versus left cosets and right multiplication by voltages in \eqref{eq:def}; the latter agrees with the way ordinary voltage graphs and lifts have been introduced in the monograph by Gross and Tucker \cite{gt87}. We refer the interested reader to \cite{pt21} for more details; here, we state a version of Theorem 4 of \cite{pt21} that applies to our setting.

\begin{theorem}
Let $\Gamma^*$ be a graph and let $G$ be a group of automorphisms of $\Gamma^*$ that acts freely on the arcs of the graph. Then, on the quotient graph $\Gamma=\Gamma^*/G$ with vertex set $V$ and arc set $A$, there exists a voltage assignment $\alpha: A\to G$ with the property that $\alpha(a^-) = \alpha(a)^{-1}$ and a function $\omega$ that assigns to every $u\in V$ a subgroup $G_u$ of $G$, such that the factored lift $\Gamma^{(\alpha, \omega)}$ is isomorphic to the original graph $\Gamma^*$.
\end{theorem}

Thus, factored lifts are, on the one hand, a generalisation of ordinary voltage lifts, dealing with a group of automorphisms acting freely on {\em vertices}, to the situation of a group of automorphisms acting freely on {\em arcs}. On the other hand, factored lifts are a special case of the general lifts of \cite{pt21} that allow for an arbitrary action of a group of automorphisms. For completeness, there is another well-known generalisation of ordinary lifts, which are the so-called {\em relative} or {\em permutation lifts}, see Gross and Tucker \cite{gt87} and Dalf\'o,  Fiol, Pavl\'ikov\'a and  \v{S}ir\'a\v{n} \cite{dfps21}. The generalisation in the latter, however, is in terms of coverings ---it extends consideration of regular coverings arising from ordinary lifts to general coverings--- and \textit{does not refer to group actions as in the present paper}.

\section{Counting lifts of walks}\label{sec:walks}

Lifts of walks are of central importance in the study of coverings. In an ordinary lift $\Gamma^{\alpha}$ of a base graph $\Gamma$ with a voltage assignment $\alpha$ in a group $G$, every walk $W$ in $\Gamma$ starting at a vertex $u$ lifts, for each $h\in G$, to a unique walk $\widetilde W_h$ in $\Gamma^{\alpha}$ starting at the vertex $(u,h)$. Here `lifting' means that the projection $\pi:\ \Gamma^{\alpha}\to \Gamma$ given by erasing the group coordinate maps arcs of the walk $\widetilde W_h$ bijectively onto those of $W$. The situation is, however, a bit different in a factored lift $\Gamma^{(\alpha,\omega)}$ arising from $\Gamma$ by a combined voltage assignment $(\alpha,\omega)$ in $G$, and this is what we now aim to explain.

By the definition of a factored lift, for every $h\in G$, an arc $a\in \overrightarrow{uv}$ of $\Gamma$ lifts to the arc $(a,h)$ in $\Gamma^{(\alpha,\omega)}$ emanating from the vertex $(u,hG_u)$ and terminating at the vertex $(v,h\alpha(a)G_v)$. But for a fixed $h$ and an arbitrary $g\in G_u$, every arc of the form $(a,hg)$ emanates from the {\em same} vertex $(u,hG_u)$ since $hgG_u=hG_u$. This way, one obtains $|G_u|$ arcs $(a,hg)$ in the factored lift, all emanating from the vertex $(u,hG_u)$ and projecting onto the arc $a$ of $\Gamma$ by the same projection $\pi$ as above.  Note, however, that their terminal vertices $(v, hg\alpha(a)G_v)$ may be different, because of the fact that for various $g\in G_u$ one may end up with {\em different} cosets $hg\alpha(a)G_v$.

\begin{remark}\label{rem:beta} A useful way to look at lifts of arcs by combined voltage assignments is to imagine that every arc $a\in \overrightarrow{uv}$ of the base graph, carrying voltage $\alpha(a)\in G$, has been assigned a {\em set} $\overline{\alpha}(a)$ of voltages, of the form $\overline{\alpha}(a) = G_u\alpha(a) = \{g\alpha(a);\ g\in G_u\}$. This set of voltages $\overline{\alpha}(a)$ may then be considered to induce the $|G_u|$ lifts $(a,hg)$ of the arc $a$ for $g\in G_u$, all emanating from the same vertex $(u,hG_u)$ but terminating at some vertex contained in the {\em set} of vertices denoted by  $(v,h \overline{\alpha} (a)G_v)$ and defined by $(v,h\overline{\alpha}(a)G_v) = \{(v,hg\alpha(a)G_v);\ g\in G_u\}$.
\end{remark}

This feature propagates when one continues to follow arcs along a walk in the base graph $\Gamma$. To see this, let $W=a_1a_2\ldots a_\ell$ be a walk of length $\ell$ in $\Gamma$ consisting of $\ell$ consecutive arcs, where $a_j\in \overrightarrow {v_jv_{j+1}}$ for $j\in [\ell] = \{1,2,\ldots, \ell\}$. Recalling the notation $\omega(u)=G_u$ for every vertex $u$ of $\Gamma$, let us choose $\ell$ elements $g_j\in \omega(v_j)$ for $j\in [\ell]$ in an arbitrary way. Further, let $h_1=1$ and for $j\in [\ell]$ we recursively define $h_{j+1} = h_jg_j \alpha(a_j)$. Then, the walk $W$ lifts to a walk $\widetilde W$ in $\Gamma^{(\alpha,\omega)}$ of the form
\begin{equation}\label{eq:walk-lift}
\widetilde{W} = (a_1,h_1g_1)(a_2,h_2g_2)\ldots (a_j,h_jg_j) \ldots (a_\ell,h_{\ell}g_{\ell}),
\end{equation}
where, for every $j\in [\ell]$, the arc $(a_j,h_jg_j)$ of $\Gamma^{(\alpha,\omega)}$ starts at the vertex $(v_j,h_j\omega(v_j))$ and ends at the vertex $(v_{j+1},h_jg_j\alpha(a_j)\omega(v_{j+1}))$ coinciding with $(v_{j+1},h_{j+1}\omega(v_{j+1}))$ due to the definition of the sequence $(h_j)_{j\in [\ell]}$. It may be checked that every lift of a walk arises this way. In particular, our walk $W$ gives rise to $\prod_{j\in [\ell]}|\omega(v_{j})|$ lifts $\widetilde W$ in $\Gamma^{(\alpha,\omega)}$ as in \eqref{eq:walk-lift}, each projecting onto $W$ by $\pi$; their count matches the number of choices of elements $g_j\in \omega(v_j)$ for $j\in [\ell]$. Note that if $\omega$ is a trivial assignment, this gives a generalisation of the unique-walk-lifting property for ordinary lifts.

In the special case when the walk $W$ in $\Gamma$ as above is {\em closed}, that is, when $v_1=v_{\ell+1}$, the lift $\widetilde W$ given by \eqref{eq:walk-lift} is a closed walk in $\Gamma^{(\alpha,\omega)}$ if and only if $\omega (v_{\ell+1})=\omega(v_1)$ and, at the same time, $h_{\ell+1}\in \omega(v_1)$. But by our recursion (with $h_1=1$), the last condition means that
\begin{equation}\label{eq:walk-recurs}
h_{\ell+1} = g_1\alpha(a_1)g_2\alpha(a_2)\ldots g_\ell\alpha(a_\ell) \in \omega(v_1)\ .
\end{equation}
Further, as $g_1$ already belongs to $\omega(v_1)$, it follows from \eqref{eq:walk-recurs} that for every $g\in \omega(v_1)$ there is an element $\bar g_1\in \omega(v_1)$ such that
\begin{equation}\label{eq:walk-recurs-e}
\bar g_1\alpha(a_1)g_2\alpha(a_2)\ldots g_\ell\alpha(a_\ell) = g \ .\end{equation}
 
Equation \eqref{eq:walk-recurs-e} may be usefully interpreted by recalling the voltage sets $\overline{\alpha}(a) = G_u\alpha(a) = \omega(u)\alpha(a) = \{g\alpha(a);\ g\in G_u\}$ defined for arcs $a$ emanating from a vertex $u$ in the base graph. Namely, if one lets $\overline{\alpha}(a_j)= \omega(v_j)\alpha(a_j)$ for $j\in [\ell]$, then the multiset $\overline{\alpha}(W) = \overline{\alpha}(a_1) \overline{\alpha}(a_2) \ldots \overline{\alpha} (a_\ell)$, consisting of all products of the form $\omega_1\omega_2\cdots\omega_\ell$ with $\omega_j\in \overline {\alpha}(a_j)$ for $j\in [\ell]$, is the multiset of the total voltages (also known as `net voltages') of $W$ which one accumulates by moving along arcs of $W$ in $\Gamma$ and multiplying the respective voltages on these arcs. The equation \eqref{eq:walk-recurs-e} then gives a one-to-one correspondence between the multiset of those products in $\overline{\alpha}(W)$ that are equal to $g\in \omega(v_1)$ and the multiset of the products in $\overline {\alpha}(W)$ that are equal to any other element of $\omega(v_1)$. It follows that there is a one-to-one correspondence between the set of closed walks in $\Gamma$ of length $\ell$ rooted at $v_1$ that have net voltage equal to $g\in G_u$ on the one hand, and those of net voltage $1$ on the other hand. We summarize this observation for future reference.

\begin{remark}\label{rem:walks}
Let $\Gamma$ be a base graph equipped with a combined voltage assignment $(\alpha,\omega)$ in a group $G$. In the factored lift $\Gamma^{(\alpha, \omega)}$, lifts of closed walks $W$ of length $\ell$ in the base graph $\Gamma$, rooted at $u$ and of net voltage $g$ for any particular $g\in G_u$, are in a one-to-one correspondence with lifts of walks $W$ with the same parameters but with net voltage $1$.
\end{remark}

Our intended study of lifts of eigenvectors and eigenvalues from a combined voltage graph to a factored lift requires introducing further notation, the origins of which are found in Dalf\'o, Fiol, and \v{S}ir\'a\v{n} \cite{dfs19}, and Dalf\'o, Fiol, Pavl\'ikov\'a, and \v{S}ir\'a\v{n} \cite{dfps21}. 
Given a combined voltage graph $\G=(V,A)$ of order $k$ under the assignment $(\alpha, \omega)$ in a voltage group $G$ we first assign to it a $k\times k$ matrix $\B = \B(\G;\alpha,\omega)$ indexed with the set $V$, entries of which are elements of the complex group algebra $\CC[G]$ of $G$. For every $u\in V$, we first introduce a specific element $G_u^+ \in \CC[G]$ by letting  $G_u^+=\sum_{g\in G_u}g$. With this in hand, for every $u,v\in V$, the $(u,v)^{\rm th}$ element of $\B$ is defined by
\begin{equation}
\label{eq:BG}
\B_{u,v} = \sum_{a\in \overrightarrow{uv}}G_u^+\alpha(a) = G_u^+\sum_{a\in \overrightarrow{uv}}\alpha(a) \ ,
\end{equation}
where, as before, $\overrightarrow{uv}$ is the set of all arcs from $u$ to $v$ in $\G$, with $\B_{u,v}=0$ if $\overrightarrow{uv}=\emptyset$.

The matrix $\B = \B(\G;\alpha,\omega)$ associated with a combined voltage graph $\Gamma$ by \eqref{eq:BG}  enables us to determine the number of closed walks of a given length and rooted at a given vertex in the factored lift $\Gamma^{(\alpha, \omega)}$ as follows.

\begin{proposition}
\label{prop:walks}
Assume that, for a given $\ell\ge 1$, the $(u,u)^{\rm th}$ entry of the $\ell^{\rm th}$ power $\B^\ell$ of the matrix $\B=\B(\G;\alpha,\omega)$ is equal to the element $(\B^{\ell})_{uu}=\sum_{g\in G}b_{g}^{(\ell)}g$ of the group algebra $\CC[G]$. Then, the number $n(u,\ell)$ of closed walks in the factored lift $\G^{(\alpha,\omega)}$,  rooted at a vertex $(u,G_u)$, is equal to
\begin{equation}
n(u,\ell)=|G_u|\cdot b_{e}^{(\ell)},
\label{eq:walks}
\end{equation}
where $e$ is the identity element of $G$.
\end{proposition}

\begin{proof}
We begin by pointing out that, by \eqref{eq:BG}, every entry of the $u^{\rm th}$ row of $\B$ is a left multiple of the element $G_u^+=\sum_{g\in G_u}g$ of the group algebra $\CC[G]$. Because of this, every arc $a$ emanating from a vertex $u$ in $\G$ may be viewed as being  equipped with a {\em set} of voltages of the form $\overline{\alpha}(a)= G_u\alpha(a) = \{g\alpha(a);\ g\in G_u\}$ from Remark \ref{rem:beta}. This set of voltages may, in turn, be identified with the element $G_u^+\alpha(a) \in \CC[G]$ constituting one term in the definition of the entry $\B_{u,v}$ for a specific arc $a\in \overrightarrow{uv}$. For the rest of the argument, assume that a vertex $u$ of $\Gamma$ has been fixed.

Invoking Remark \ref{rem:beta} again and making use of the obvious interpretation of entries of a power of a matrix, it follows that for $g\in G$ the coefficient $b_{g}^{(\ell)}$ of the group-algebra element $(\B^{\ell})_{uu}= \sum_{g\in G}b_{g}^{(\ell)}g$ is equal to the number of closed walks $W$ in $\Gamma$ of length $\ell$, rooted at $u$ and with value $g$ in the multiset $\overline{\alpha}(W)$. In particular, for $g$ equal to the unit element $1\in G$, the coefficient $b_{1}^{(\ell)}$ counts the number of closed walks $W$ in $\Gamma$ of length $\ell$, rooted at $u$ but with trivial net voltage $1\in \overline{\alpha}(W)$. By Remark \ref{rem:walks}, however, there is a one-to-one correspondence between closed walks of length $\ell$ rooted at $(u,G_u)$ in the factored lift, with projections onto the base graph having net voltages $1$ and $g$, respectively, for an arbitrary $g\in G_u$. This translates to the fact that $b_{1}^{(\ell)} = b_{g}^{(\ell)}$ for every $g\in G_u$ and so the number of closed walks of length $\ell$ in the factored lift, rooted at a vertex $(u,G_u)$, is equal to $|G_u|\cdot b_1^{(\ell)}$, as claimed.  \end{proof}

As an illustration of Proposition \ref{prop:walks}, consider again the example of Figure \ref{fig1}, but replacing the additive representation of the cyclic voltage group of order $4$ by its multiplicative version $G= \langle g\ |\ g^4=1\rangle$, with $G_u=\{1\}$ and $G_v=\{1,g^2\}$. The matrix $\B$ from \eqref{eq:BG} associated with the combined voltage graph on the left-hand side of Figure \ref{fig1} has the form
\begin{equation}
\B=\left(
\begin{array}{cc}
g+g^{-1} & 1+g\\
(1+g^2)(1+g^{-1}) & 0
\end{array}
\right)=\left(
\begin{array}{cc}
g+g^{-1} & 1+g\\
1+g+g^{-1}+g^2& 0
\end{array}
\right).
\label{B-superoctahedron}
\end{equation}
Taking $\ell=5$, one may check that $(\B^5)_{uu}=176g^3+160g^2+176g+160$ and $(\B^5)_{vv}=80g^3+80g^2+80g+80$. By Proposition \ref{prop:walks}, for the number of closed walks rooted at $(u,G_u)$ and at $(v,G_v)$ in the factored lift one obtains $n(u,5)=|G_u|\cdot b_1^{(5)}=1\cdot160=160$ and $n(v,5)=|G_v|\cdot b_1^{(5)}=2\cdot 80=160$. The two values coincide as the factored lift is vertex-transitive.


\section{Voltage group representations, lifts of spectra}
\label{sec:factored-lifts}

To explain connections between factored lifts and representations of voltage groups, let $\Gamma=(V,A)$ be a base graph of order $k$ with a combined voltage assignment $(\alpha,\omega)$ in a group $G$. Let $\rho$ be a complex irreducible representation of $G$ in $\mathbb{C}^d$ of dimension $d=d(\rho)$. Recalling the matrix $\B=\B(\G; \alpha, \omega)$ defined by \eqref{eq:BG} in the previous section, we now link this matrix with the representation $\rho$ by introducing a $dk\times dk$ complex block matrix $\B(\rho)$. For every ordered pair $(u,v)$ of vertices of $V$ the  $(u,v)^{\rm th}$ block entry of $\B(\rho)$ is defined to be the $d\times d$ matrix
\begin{equation}
\label{eq:Brho}
\B_{u,v}(\rho) = \sum_{a\in \overrightarrow{uv}}\ \sum_{h\in G_u}\ \rho(h\alpha(a)),
\end{equation}
where the sum of matrices is defined in the usual way; the block entry $\B_{u,v}(\rho)$ is the all-zero $d\times d$ matrix if $\overrightarrow{uv}=\emptyset$. (We assume throughout that the indexation within the $d\times d$ blocks of $\B_{u,v}(\rho)$ by the set $\{1,2,\ldots,d\}$ is the same across all the $k^2$ blocks of this kind, which themselves are indexed by pairs of elements of $V$.)

In order to simplify the forthcoming calculations, for $G_u$ and its $\CC[G]$-variant $G_u^+ = \sum _{g\in G_u} g$ introduced in the previous section, we let
\begin{equation}\label{eq:Brhosimple}
\rho(G_u) = \rho(G_u^+) = \sum_{h\in G_u}\rho(h) \ .
\end{equation}
Combined with the fact that $\rho$ is a group homomorphism, the notation of \eqref{eq:Brhosimple} enables one to rewrite the defining equation \eqref{eq:Brho} in the form
\begin{equation}
\label{eq:Brho1}
\B_{u,v}(\rho) = \rho(G_u) \sum_{a\in \overrightarrow{uv}} \rho(\alpha(a)),
\end{equation}
which can be advantageously interpreted by saying that the $u^{\rm th}$ row of the matrix $\B(\rho)$ is a left multiple of the $d\times d$ factor $\rho(G_u)$ which is `constant' for any fixed $u\in V$.

In fact, when $\rho_0$ is the trivial representation, then $\B(\rho_0)$ is a quotient matrix of a regular (or equitable) partition of the factored lift graph, where cells correspond to fibres $F(u)=\{(u,H)\in V^{(\alpha, \omega)}\ |\ H\in G/G_u\}$ introduced in the previous section; see the example in Figure \ref{fig:factored-lift}. In particular, the largest eigenvalue of $\B(\rho_0)$ corresponds to the spectral radius of the factored lift.

To work with group representations, let $\Irep(G)$ be a complete set of irreducible representations of a finite group $G$. Let $H<G$ be an arbitrary subgroup of $G$ and, for any $\rho \in \Irep(G)$ with dimension $d(\rho)$, let $\rho(H) = \sum_{h\in H} \rho(h)$. That is, $\rho(H)$ is the sum of $d(\rho)$-dimensional complex matrices $\rho(h)$ taken over all elements $h\in H$. We note that, in general, $\rho(H)$ may be the zero-matrix, although, of course, all the matrices $\rho(h)$ for $h\in H$ are non-singular. We make use of the following result of \cite[Prop. 3.3]{dfps21} obtained earlier by the authors of the present paper.

\begin{proposition}[\cite{dfps21}]
\label{prop:ours}
For every group $G$ and every subgroup $H<G$ of index $n=[G:H]$ one has
$$
\sum_{\rho\in \Irep(G)}
d(\rho)\cdot\rank(\rho(H)) = n.
$$
\end{proposition}

We also need some preparation for working with column vectors of dimension a multiple of $d$, say, $d\ell$ for some $\ell\ge 1$. A complex vector $\f$ of dimension $d\ell$ is represented in the form $\f^{\top}=(\f_1,\f_2, \ldots,\f_\ell)^{\top}$, where, for $t\in [\ell] = \{1,2,\ldots,\ell\}$, each $\f_t$ is a $d$-dimensional column vector called a {\em $d$-segment} of $\f$. For such a vector $\f$ of dimension $d\ell$ and for every $j\in [d]=\{1,2,\ldots, d\}$, the {\em $j$-section} of $\f$ is the $\ell$-dimensional column vector $\f_{[j]}$ given by $\f_{[j]}^{\top}= (\f_{1,j}, \f_{2,j}, \ldots,\f_{\ell,j})^{\top}$, where $\f_{t,j}$ is the $j$-th coordinate of $\f_t$ for each $t\in [\ell]$.

Suppose now that a $dk$-dimensional complex column vector $\f$ is an eigenvector of our $dk\times dk$ matrix $\B(\rho)$ for some complex eigenvalue $\lambda$. Assuming consistent indexation of $d\times d$ blocks of $\B(\rho)$ and $d$-segments of $\f$ by the vertex set $V$, for every $v\in V$ the $v^{\rm th}$ $d$-segment of $\f$ is denoted $\f(v)$.

It remains to introduce a condition for eigenvectors, which we refer to in what follows. We say that, for a $d$-dimensional representation $\rho$ of $G$, an eigenvector $\f$ for some eigenvalue of the $dk\times dk$ dimensional matrix $\B(\rho)$ satisfies the condition ${\bf (C)}$ if the following is fulfilled:
\begin{equation}
\label{eq:Cond-d}
{\bf (C)}\quad \mbox{For every   $v\in V$,  $\f(v)\ne \vec0$ implies $\rho(h)\f(v) = \c_v$ for every  $h\in G_v$,}
\end{equation}
where $\c_v$ is a $d$-dimensional column vector that may depend on $v$ but is constant over the elements $h\in G_v$. The important fact is that the condition ${\bf (C)}$ implies that the $dk^\omega$-dimensional vector $\f^+$, with $k^\omega$ $d$-segments indexed by $k^\omega$ ordered pairs $(v,hG_v)$ for $v\in V$ and $hG_v\in G/G_v$, given by $\f^+(v,hG_v) = \rho(h)\f(v)$, is well defined and does not depend on representatives of cosets in $G/G_v$.

We are now ready to state and prove our first result, linking eigenvectors and eigenvalues of a factored lift with representations of its voltage group.

\begin{theorem}
\label{theo-sp-d-dim}
Let $\G$ be a graph on a set $V$ of $k$ vertices, with a combined voltage assignment $(\alpha,\omega)$ in a finite group $G$ that assigns a subgroup $G_v$ of $G$ to every vertex $v\in V$, and let $k^\omega = \sum_{v\in V}[G:G_v]$ be the number of vertices of the factored lift $\G^{(\alpha,\omega)}$. Let $\rho$ be a complex irreducible representation of $G$ of dimension $d\ge 1$ and let $\B(\rho)$ be the associated complex $dk\times dk$ matrix defined by \eqref{eq:Brho}. Further, let $\f$ be a $dk$-dimensional eigenvector of $\B(\rho)$ for some complex eigenvalue $\lambda$ of $\B(\rho)$, with $d$-blocks $\f(v)$ for $v\in V$, which fulfils the condition {\rm {\bf (C)}}. Then:
\begin{itemize}
\item[$(i)$]
The $dk^\omega$-dimensional vector $\f^+$ with $k^\omega$ $d$-segments indexed by $k^\omega$ ordered pairs $(v,hG_v)$ for $v\in V$ and $hG_v\in G/G_v$, given by $\f^+(v,hG_v) = \rho(h)\f(v)$, is well defined and does not depend on representatives of cosets in $G/G_v$.
\item[$(ii)$]
For every $j\in [d]$, the $j$-section $\f^+_{[j]}$ of $\f^+$ is a $k^\omega$-dimensional eigenvector of the factored lift $\G^{(\alpha,\omega)}$ for the same eigenvalue $\lambda$ as above.
\item[$(iii)$]
Let $S$ be the system of $k^\omega$ linear equations of the form ${\bf c}\rho(h)\f(v) = 0$ for $v\in V$ and $hG_v\in G/G_v$ for an unknown row vector ${\bf c}=(c_1,c_2,\ldots,c_d)$ of dimension $d$. The set of $j$-sections $\{\f^+_{[j]}\ |\ j\in [d]\}$ is linearly independent if and only if the system $S$ has only a trivial solution ${\bf c} = {\bf 0}$. In particular, this is satisfied if $V$ contains a subset $U$ of $d$ vertices such that the $d$-segments $\f(u)$ for $u\in U$ are linearly independent.
\end{itemize}
\end{theorem}

\begin{proof}
Part $(i)$ is a consequence of the condition {\bf (C)}, and so we move on to $(ii)$. Our assumption that $\f$ is a column eigenvector of $\B(\rho)$ for an eigenvalue $\lambda$ is, with the help of \eqref{eq:Brho1}, equivalent to stating that the $d$-segments $\f(v)$ of $\f$ for $v\in V$ satisfy
\begin{equation}
\label{eq:Bfd}
\lambda \f(u) = \sum_{v\sim u}\B(\rho)_{u,v}\f(v) = \sum_{v\sim u}\ \rho(G_u) \sum_{a\in \overrightarrow{uv}} \rho(\alpha(a)) \f(v)
\end{equation}
for every vertex $u\in V$. We also assume that $\f$ satisfies our assumption {\bf (C)} as stated in \eqref{eq:Cond-d}. For such a $\lambda$ and $\f$, we introduce a new column vector $\f^+$ of dimension $dk^\omega$ whose $d$-segments, indexed by the $k^\omega$ pairs $(v,gG_v)\in V^{(\alpha,\omega)}$, are defined by
\begin{equation}\label{eq:Fd}
\f^+(v,hG_v)= \rho(h)\f(v) {\rm\ for\ every\ } v\in V {\rm \ and\ every\ } h\in G\ ;
\end{equation}
the fact that $\f^+$ is well defined is a direct consequence of the assumption {\bf (C)}. Multiplying \eqref{eq:Bfd} by the $d\times d$ matrix $\rho(g)$ from the left and using \eqref{eq:Brho} with \eqref{eq:Brho1} gives
\begin{equation}
\label{eq:Bfgd}
\lambda \rho(g)\f(u) =  \sum_{v\sim u}\ \sum_{a\in \overrightarrow{uv}}\ \sum_{h\in gG_u}\ \rho(h)\rho(\alpha(a)) \f(v)\ .
\end{equation}
Using now the equations \eqref{eq:Fd}, one sees that \eqref{eq:Bfgd} is equivalent to the statement that
\begin{equation}
\label{eq:Bf+gd}
\lambda\f^+(u,gG_u) = \sum_{v\sim u}\sum_{a\in \overrightarrow{uv}}\sum_{h\in gG_u}\f^+(v,h\alpha(a)G_v)\ .
\end{equation}
The important conclusion that follows from \eqref{eq:Bf+gd} is that, for every $j\in [d]$, the $j$-section $\f^+_{[j]}$ of $\f$ is a $k^\omega$-dimensional vector satisfying the equation
\begin{equation}
\label{eq:Bf+gd2}
\lambda\f^+_{[j]}(u,gG_u) = \sum_{v\sim u}\sum_{a\in \overrightarrow{uv}}\sum_{h\in gG_u}\f^+_{[j]}(v,h\alpha(a)G_v)\ . 
\end{equation}
But \eqref{eq:Bf+gd2} demonstrates that, for every $j\in [d]$, the $j$-section $\f^+_{[j]}$ is an eigenvector of the factored lift $\G^{(\alpha,\omega)}$ that belongs to the {\em same} eigenvalue $\lambda$ we started with.
\smallskip

For part $(iii)$, the $j$-sections $\f^+_{[j]}$ for $j\in [d]$ form a linearly independent set if and only if the linear combination $c_1\f^+_{[1]} + c_2\f^+_{[2]} + \cdots + c_d\f^+_{[d]}$ results in a $k^\omega$-dimensional zero vector only in the trivial case, that is, when ${\bf c} = (c_1,c_2,\ldots,c_d)$ is a zero row vector. This must hold for every $d$-segment of the $j$-sections, but since by \eqref{eq:Fd} the $(u,gG_u)$-th coordinate of $\f^+_{[j]}$ is the $j$-th coordinate $(\rho(g)\f(u))_j$ of the $d$-segment $\rho(g)\f(u)$, the above linear combination equates to a zero vector if and only if the following system of $k^\omega$ equations (for every vertex $(u,gG_u)$ of the factored lift)
\begin{equation}
\label{eq:lincomb}
c_1(\rho(g)\f(u))_1 + c_2(\rho(g)\f(u))_2 + \cdots + c_d(\rho(g)\f(u))_d = 0
\end{equation}
has only the trivial solution, namely, the zero row vector ${\bf c}$. But the system \eqref{eq:lincomb} may simply be rewritten in the form ${\bf c}{\cdot}\rho(g){\cdot}\f(u) = 0$. This implies the validity of the statement $(iii)$, including the particular observation about a subset $U$ of $d$ linearly independent $d$-segments $\f(u)$ for $u\in U$.
\end{proof}

If one is interested only in calculating the spectrum $\spec(\G^{(\alpha,\omega)})$ of a factored lift, it turns out that it is sufficient to consider the spectra $\spec(\B(\rho))$ of the complex matrices $\B(\rho)$, taken over a complete set of irreducible complex representations $\rho$ of the voltage group $G$, as our next result shows.

\begin{theorem}
\label{th:spectd}
Let $(\alpha,\omega)$ be a combined voltage assignment on a graph $\G=(V,A)$ with $k$ vertices in a group $G$ and let $n_u=|G:G_u|$ for every $u\in V$, and with the order of the factored lift $\G^{(\alpha,\omega)}$ equal to $k^{\omega}= \sum_{u\in V}n_u$. Let $G$ have order $n$, with $\nu$ conjugacy classes, and let $\{\rho_r : \ r=0,1,\ldots,\nu-1\}$ be a complete set of complex irreducible representations of $G$, of dimensions $d(\rho_r)=d_r$, so that $\sum_{r=0}^{\nu-1}d_r^2=n$. Let ${\cal B}$ be the multiset of eigenvalues
$$
{\cal B} = \bigcup_{r=0}^{\nu-1}d_r\,\spec(\B(\rho_r)).
$$
of cardinality $\sum_{r=0}^{\nu-1}kd_r^2=kn$. Then, the following statements hold:
\begin{itemize}
\item[$(i)$]
The multiset ${\cal B}$ contains at most $k^{\omega}=\sum_{u\in V}n_u$ non-zero eigenvalues.
\item[$(ii)$]
The spectrum of the factored lift $\G^{(\alpha,\omega)}$ is the multiset ${\cal B}\setminus {\cal Z}$, where ${\cal Z}$ is a multiset containing $kn-k^{\omega}$ zeros.
\end{itemize}
\end{theorem}

\begin{proof}
For part $(i)$, let $u\in V$ be a vertex with voltage subgroup $\omega(u)=G_u<G$. Then, for a given irreducible representation $\rho_r$ of dimension $d_r$, the $u^{\rm th}$ block-row of $\B(\rho_r)$, which is a $d_r\times d_rk$ matrix denoted by $\B(\rho_r)_u$, satisfies, by equation~\eqref{eq:Brho1},
$$
\rank\bigl(\B(\rho_r)_u\bigr)
\leq
\rank\bigl(\rho_r(G_u)\bigr),
$$
where
$$
\rho_r(G_u)=\sum_{h\in G_u}\rho_r(h).
$$
Since $\B(\rho_r)$ is obtained by stacking the block-rows $\B(\rho_r)_u$, for $u\in V$, we have
$$
\rank\bigl(\B(\rho_r)\bigr)
\leq
\sum_{u\in V}\rank\bigl(\B(\rho_r)_u\bigr)
\leq
\sum_{u\in V}\rank\bigl(\rho_r(G_u)\bigr).
$$
Moreover, since $\B(\rho_r)$ is a square $d_rk\times d_rk$ matrix, the number of its non-zero eigenvalues, counted with algebraic multiplicity, is at most $\rank(\B(\rho_r))$. Therefore, the number of non-zero eigenvalues contributed to ${\cal B}$ by $\B(\rho_r)$, taking into account its multiplicity $d_r$ in ${\cal B}$, is at most
$$
d_r\sum_{u\in V}\rank\bigl(\rho_r(G_u)\bigr).
$$
Summing over all irreducible representations and using Proposition~\ref{prop:ours}, we obtain
$$
\sum_{r=0}^{\nu-1}
d_r\sum_{u\in V}\rank\bigl(\rho_r(G_u)\bigr)
=
\sum_{u\in V}\sum_{r=0}^{\nu-1}
d_r\rank\bigl(\rho_r(G_u)\bigr)\\
=
\sum_{u\in V}[G:G_u]\\
=
k^{\omega}.
$$
Thus, the multiset ${\cal B}$ contains at most $k^{\omega}$ non-zero eigenvalues, as claimed.

To establish $(ii)$, we use the irreducible characters $\chi_r$ associated with each irreducible representation of $G$. Let $\A$ be the adjacency matrix of $\G^{(\alpha,\omega)}$, with entries $a_{(u,G_u)(v,G_v)}$ for $u,v\in V$. It is well known that the total number of rooted closed walks of length $\ell$ in the factored lift is equal to the trace of the $\ell^{\rm th}$ power $\A^\ell$ of $\A$, with elements $a^{(\ell)}_{(u,G_u)(v,G_v)}$. But the same trace is also equal to the sum of the $\ell^{\rm th}$ powers of the eigenvalues (including multiplicities) of $\A^\ell$. This implies that
\begin{equation}\label{eq:trace}
\tr(\A^{\ell}) =\sum_{u\in V}n_u a_{(u,G_u)(u,G_u)}^{(\ell)} =\sum_{\lambda\in \spec(\G^{(\alpha,\omega)})}\lambda^{\ell},
\end{equation}
where the middle part of \eqref{eq:trace} is a consequence of the fact that left multiplication by elements of $G$ induces automorphisms of the factored lift that act transitively on its fibres.

By Proposition \ref{prop:walks}, the $u^{\rm th}$ diagonal entry of $\A^\ell$ is a $|G_u|$-multiple of the coefficient $b_e^{(\ell)}$ at the `identity term' $e$ of the $\CC[G]$-group ring element $(\B^{\ell})_{uu} = \sum_{g\in G}b_g^{(\ell)}g$. On the other hand, in \cite{dfps21} (by the same authors of the present paper), it was proved that, in terms of characters, the same
`identity' coefficient $b_e^{(\ell)}$ admits the following evaluation:
\begin{equation}\label{eq:char}
b_e^{(\ell)}=\frac{1}{n}\sum_{r=0}^{\nu-1}d_r\chi_r((\B^{\ell})_{uu}).
\end{equation}
Combining \eqref{eq:char} with the aforementioned result of Proposition \ref{prop:walks} gives
$$
a_{(u,G_u)(u,G_u)}^{(\ell)}=|G_u|b_e^{(\ell)}=\frac{1}{n_u}\sum_{r=0}^{\nu-1}d_r\chi_r((\B^{\ell})_{uu}).
$$
Consequently, putting the pieces together and using $\mu$ as a variable for eigenvalues of the matrices $\B(\rho_r)$, one obtains
\begin{align*}
\sum_{\lambda\in \spec(\G^{(\alpha,\omega)})}\lambda^{\ell}& =\tr(\A^{\ell})
= \sum_{u\in V}n_u a_{(u,G_u)(u,G_u)}^{(\ell)}
=\sum_{u\in V}
\sum_{r=0}^{\nu-1} d_r\chi_r((\B^{\ell})_{uu}) \\
 &=\sum_{r=0}^{\nu-1}d_r\sum_{u\in V}\tr(\B(\rho_r)^{\ell})
  =\sum_{r=0}^{\nu-1}d_r\sum_{\mu\in\spec(B(\rho_r)^{\ell})}\mu^{\ell},
\end{align*}
where the leftmost sum contains $k^{\omega}=\sum_{u\in V}n_u$ terms, as opposed to the formally $nk$ terms of the rightmost sum, but by part $(i)$ one is free to remove $kn-k^{\omega}$ zeros from the last sum without affecting the equality. Hence, since such `adjusted equalities' hold for every $\ell=1,\ldots,k^{\omega}-1$, the multisets of eigenvalues of the factored lift $\G^{(\alpha,\omega)}$ and those in ${\cal B}$ must coincide up to a multiset ${\cal Z}$ of $kn-k^\omega$ zeros (see, for instance, Gould \cite{g99}). This completes the proof.
\end{proof}


\section{Illustrative examples}\label{sec:examples}

In this section, we give two examples illustrating issues one may encounter when considering spectra and eigenspaces of factored lifts. Since we are working with voltage groups together with their complex representations and irreducible characters, to avoid confusion, we use the symbol $e$ for the unit element of a voltage group in what follows.

\begin{figure}[ht]
\centering
\includegraphics[width=15cm]{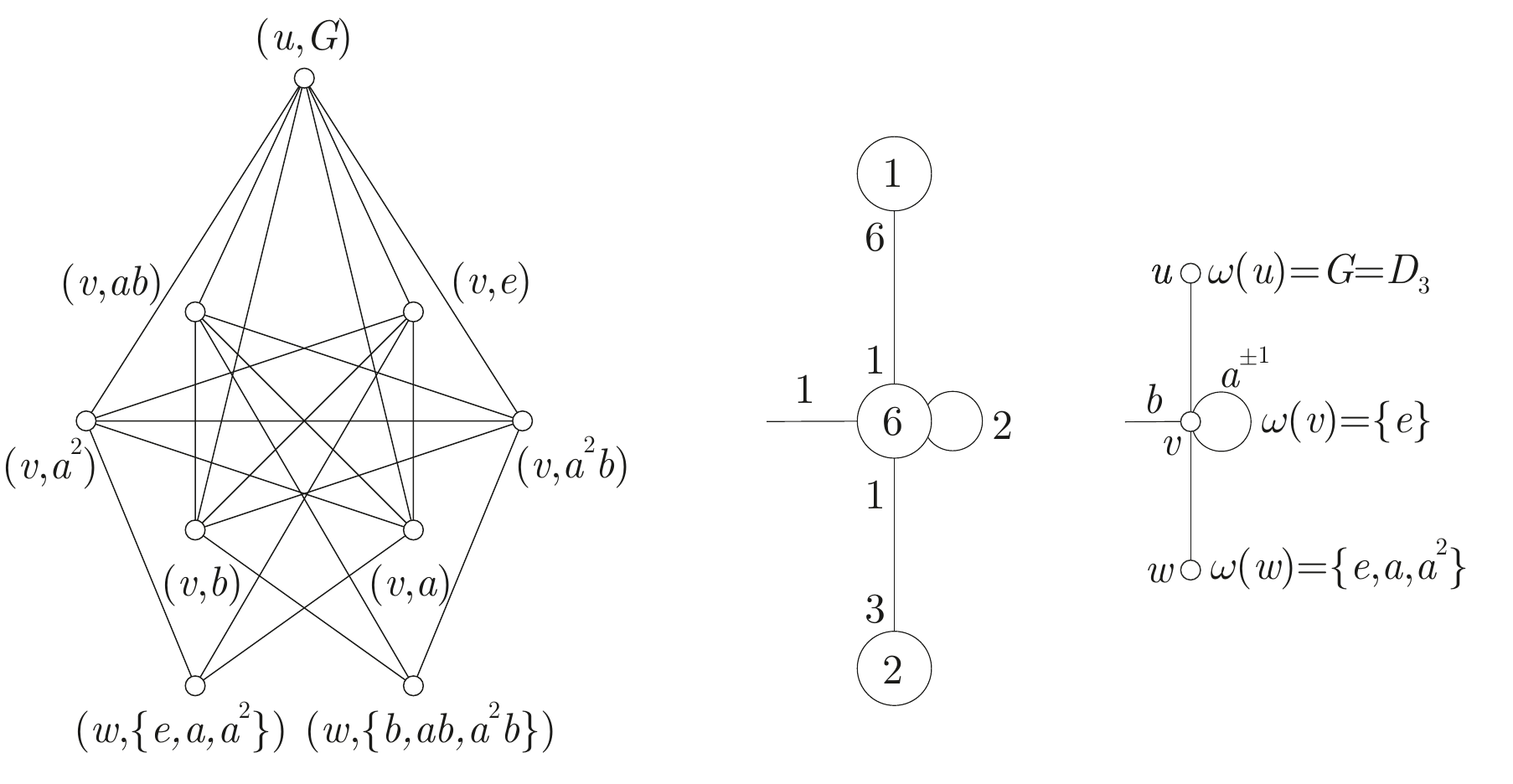}
\caption{The factored lift $\G^{(\alpha,\omega)}$ on the dihedral group $D_3$ of Example \ref{ex:example1}, a regular partition, and the corresponding combined voltage graph.}
\label{fig:factored-lift}
\end{figure}

\begin{example}
\label{ex:example1}
{\rm Consider the graph $\G^{(\alpha,\omega)}$ shown on the left-hand side of Figure \ref{fig:factored-lift}.}
\end{example}
The partition shown in the centre of Figure \ref{fig:factored-lift} indicates that the graph is a factored lift induced by an action of a dihedral group $D_3=\langle a,b\ |\ a^3= b^2 =(ab)^2=e\rangle$ of order $6$, and the corresponding combined voltage graph $\Gamma$ is shown on the right-hand side of the same figure. The matrix $\B=\B(\G)$, introduced in general by \eqref{eq:BG}, is in this case the $3\times 3$ matrix shown in \eqref{eq:dihedr-B}:
\begin{equation}\label{eq:dihedr-B}
\B=\left(
\begin{array}{ccc}
0 & \sum_{g\in D_3} g & 0\\
e & a+a^{-1}+b & e\\
0 & e + a + a^2 & 0
\end{array}
\right)
\end{equation}

A complete set of irreducible complex representations of the group $D_3$ presented above consists of the trivial and alternating $1$-dimensional representations $\rho_0$ and $\rho_1$ together with a $2$-dimensional representation $\rho_2$ with $\zeta$ a primitive complex $3$-rd root of $1$; they are all displayed in Table \ref{tab:irred-repreD3}. To save space, we use the symbols I, ${\rm Dia}(x,y)$ and ${\rm Off}(x,y)$ for the identity matrix, a diagonal matrix with entries $x,y$ (from top left to bottom right), and an off-diagonal matrix with entries $x,y$ (from top right to bottom left), all of dimension $2$.
\begin{table}[!ht]
\begin{center}
\begin{tabular}{|c|cccccc| }
\hline
& $e$  & $a$  & $a^2$  & $b$ & $ab$ & $a^2b$   \\
\hline\hline
$\rho_0$ &  $1$ &  $1$ & $1$ & $1$  & $1$ & $1$\\
\hline
$\rho_1$ &  $1$ &  $1$ & $1$ & $-1$  & $-1$ & $-1$\\
\hline
$\rho_2$ & \ \ \ I \ \ \ & ${\rm Dia}(\zeta,\zeta^{-1})$ & ${\rm Dia}(\zeta^2,\zeta^{-2})$ &
${\rm Off}(1,1)$ & ${\rm Off}(\zeta,\zeta^{-1})$ & ${\rm Off}(\zeta^2,\zeta^{-2})$\\
\hline
\end{tabular}
\caption{Irreducible representations of $D_3=\langle a,b\,|\,a^3=b^2=(ab)^2=e\rangle$.}
\label{tab:irred-repreD3}
\end{center}
\end{table}
In this example one has $\omega(u)=G_u=D_3$, $\omega(v)=G_v=\langle e\rangle$ and $\omega(w)=G_w=\langle a\rangle$, so  $\rho_0(G_u)=6$ and $\rho_1(G_u)=0$, $\rho_0(G_v)=\rho_1(G_v)=1$, and $\rho_0(G_w)=\rho_1(G_w)=3$, while, for example, $\B(\rho_1)_{u,v}=0$ and $\B(\rho_1)_{v,v} = \rho_1(G_v)(\rho_1(a)+\rho_1(a^{-1})+\rho_1(b))=1$. The $3\times 3$ matrices $\B(\rho_0)$ and $\B(\rho_1)$ are thus given by
\begin{equation}
 \B(\rho_0)=\left(
 \begin{array}{ccc}
 0 & 6 & 0\\
 1 & 3 & 1\\
 0 & 3 & 0
 \end{array}
 \right),\qquad  \B(\rho_1)=\left(
 \begin{array}{ccc}
 0 & 0 & 0\\
 1 & 1 & 1\\
 0 & 3 & 0
 \end{array}
 \right). 
 \label{eq:rho01(B)}
\end{equation}
Entries of $\B(\rho_2)$ are determined similarly, bearing in mind that this time they are $2\times 2$ blocks; for example, $\B(\rho_2)_{u,v}$ is the sum of all the six matrices appearing in the last row of Table \ref{tab:irred-repreD3}. An evaluation of $\B(\rho_2)_{u,v}$ gives
\begin{equation}
 \B(\rho_2)=\left(
 \begin{array}{cc|cc|cc}
 0 & 0 & \sum_{i=0}^2 \zeta^i & \sum_{i=0}^2 \zeta^i & 0 & 0\\
 0 & 0 & \sum_{i=0}^2 \zeta^{-i} & \sum_{i=0}^2 \zeta^{-i} & 0 & 0\\
 \hline
 1 & 0 & \zeta+\zeta^{-1} & 1 & 1 & 0\\
 0 & 1 & 1 & \zeta+\zeta^{-1} & 0 & 1\\
 \hline
 0 & 0 & \sum_{i=0}^2 \zeta^i & 0 & 0 & 0\\
 0 & 0 & 0 & \sum_{i=0}^2 \zeta^{-i} & 0 & 0\\
 \end{array}
 \right).
 \label{eq:rho2(B)}
\end{equation}

\noindent Eigenvalues and eigenvectors of $\B(\rho_r)$ for $r=0,1,2$ are listed in Table \ref{tab:eigenvv-rho012}, taking into account that the non-zero eigenvalues of $\B(\rho_0)$ and $\B(\rho_1)$ are, respectively, $\mu^{(0)}_{1,2}= 3(1\pm \sqrt{5})/2$ and $\mu^{(1)}_{1,2}= (1\pm \sqrt{13})/2$.

\begin{table}[ht]
\begin{center}
\begin{tabular}{|ccc|ccc|}
\hline
 & $\B(\rho_0)$ &  &  & $\B(\rho_1)$ & \\
\hline
\hline
 $\mu^{(0)}_1$ & $0$  & $\mu^{(0)}_2$  & $\mu^{(1)}_1$ & $0^*$ & $\mu^{(1)}_2$   \\
\hline
  $\left(
\begin{array}{c}
2 \\
\!\!\! \mu^{(0)}_1/3 \!\!\!\\
1
\end{array}\right)$ &
$\left(
\begin{array}{c}
-1\\
0\\
1
\end{array}\right)$ &
$\left(
\begin{array}{c}
2\\
\!\!\! \mu^{(0)}_2/3 \!\!\!\\
1
\end{array}\right)$ &
$\left(
\begin{array}{c}
0 \\
\!\!\! \mu^{(1)}_1/3 \!\!\!\\
1
\end{array}\right)$ & $\left(
\begin{array}{c}
-1\\
 0\\
 1
\end{array}\right)$ & $\left(
\begin{array}{c}
0\\
\!\!\! \mu^{(1)}_2/3 \!\!\!\\
1
\end{array}\right)$\\
\hline
\end{tabular}

\vskip 15pt

\begin{tabular}{|cccccc|}
\hline
 &  &   $\B(\rho_2)$   &   &  & \\
\hline
\hline
  $0$  & $0^*$  & $0^*$  & $0^*$ & $0^*$ & $-2$   \\
\hline
$\left(
\begin{array}{c}
0 \\
0\\
1\\
1\\
0\\
0
\end{array}\right)$ &
$\left(
\begin{array}{c}
0\\
-1 \\
0\\
0\\
0\\
1\\
\end{array}\right)$ &
$\left(
\begin{array}{c}
-1\\
0 \\
0\\
0\\
1\\
0\\
\end{array}\right)$ &
$\left(
\begin{array}{c}
-1\\
1 \\
0\\
1\\
0\\
0\\
\end{array}\right)$ &
$\left(
\begin{array}{c}
1\\
-1 \\
1\\
0\\
0\\
0\\
\end{array}\right)$ &
$\left(
\begin{array}{c}
0\\
0 \\
-1\\
1\\
0\\
0\\
\end{array}
\right)$\\
\hline
\end{tabular}

\end{center}
\caption{Eigenvalues and eigenvectors of the matrices $\B(\rho_r)$ for $r=0,1,2$.}
\label{tab:eigenvv-rho012}
\end{table}

Observe that for $r=0,1$, except for the eigenvalues $0$ marked by an additional star $(0^*)$, all the remaining  eigenvalues of $\B(\rho_0)$ and $\B(\rho_1)$ satisfy the condition {\bf (C)} in \eqref{eq:Cond-d} for $d=1$, since $\rho_0|_{G_u}=\rho_0|_{G_w}=\rho_1|_{G_w}=1$. This, however, does not hold for the eigenvector of $\B(\rho_1)$ that belongs to the starred eigenvalue $0^*$ because $\f(u)\neq 0$ and $\rho_1|_{G_u}\neq1$ (that is, $\rho_1(h)\f(u)$ is not constant for $h\in G_u$).

A similar phenomenon occurs for the eigenvectors of $\B(\rho_2)$ with the marked zero eigenvalues $(0^*)$, as neither  $\rho_2|_{G_u}$ nor  $\rho_2|_{G_w}$ are a non-zero constant. On the other hand, the first and last eigenvectors of $\B(\rho_2)$, that is, those corresponding to the eigenvalues $0$ and $-2$, do satisfy the condition since, in both cases, $\f(u)=\f(w)=(0,0)^{\top}$. By Theorem \ref{theo-sp-d-dim}, these eigenvalues have multiplicity $2$ in the factored lift graph.

To see in more detail what happens in the case of the eigenvalue $-2$ for $\B(\rho_2)$, the $2$-dimensional vectors constituting the corresponding $6$-dimensional eigenvector are $\f(u)=(0,0)^{\top}$, $\f(v)=(-1,1)^{\top}$, and $\f(w)=(0,0)^{\top}$; here we use transposes to save space. Following Theorem \ref{theo-sp-d-dim}, to construct the eigenvectors $\f_1^+$ and $\f_2^+$ of the factored lift, one subsequently calculates
\begin{align*}
\f^+(u,G_u) &= \f^+(w,G_w) = (0,0)^{\top}\ ,\\
\f^+(v,eG_v) &= {\rm I}(-1,1)^{\top} = (-1,1)^{\top}\ ,\\
\f^+(v,aG_v) &= {\rm Dia}(\zeta,\zeta^{-1})(-1,1)^{\top} = (-\zeta,\zeta^{-1})^{\top}\ ,\\
\f^+(v,a^2G_v)&= {\rm Dia}(\zeta^2,\zeta^{-2})(-1,1)^{\top} = (-\zeta^2,\zeta^{-2})^{\top}\ ,\\
\f^+(v,bG_v) &= {\rm Off}(1,1)(-1,1)^{\top} = (1,-1)^{\top}\ ,\\
\f^+(v,abG_v) &= {\rm Off}(\zeta,\zeta^{-1})(-1,1)^{\top} = (\zeta,-\zeta^{-1})^{\top}\ ,\\
\f^+(v,a^2bG_v) &= {\rm Off}(\zeta^2,\zeta^{-2})(-1,1)^{\top} = (\zeta^2,-\zeta^{-2})^{\top}\ .
\end{align*}
Then, taking the first or the second entries, we obtain $\f_1^{+}$ and $\f_2^{+}$, respectively:
\begin{align*}
\f_1^+&=(0,-1,-\zeta,-\zeta^2,1,\zeta,\zeta^2,0,0)^{\top};\\
\f_2^+&=(0,1,\zeta^{-1},\zeta^{-2},-1,-\zeta^{-1},-\zeta^{-2},0,0)^{\top},
\end{align*}
which can be checked to be eigenvectors for the eigenvalue $-2$ of the factored lift $\G^{(\alpha,\omega)}$ of Figure \ref{fig:factored-lift}. In summary, the spectrum of the factored lift $\G^{(\alpha,\omega)}$ of Figure \ref{fig:factored-lift} is
$$
\spec(\G^{(\alpha,\omega)})=\{\,3(1+\sqrt{5})/2,\,(1+\sqrt{13})/2,\,0^{[3]},\,(1-\sqrt{13})/2,\,3(1-\sqrt{5})/2,\, -2^{[2]}\,\}\ .
$$

\begin{example}
\label{ex:exemple3}
{\rm   Consider the graph $\G^{(\alpha,\omega)}$ shown on the left-hand side of Figure \ref{fig:factored-lift(c)}.}
\begin{figure}[ht]
    \centering
    \includegraphics[width=14cm]{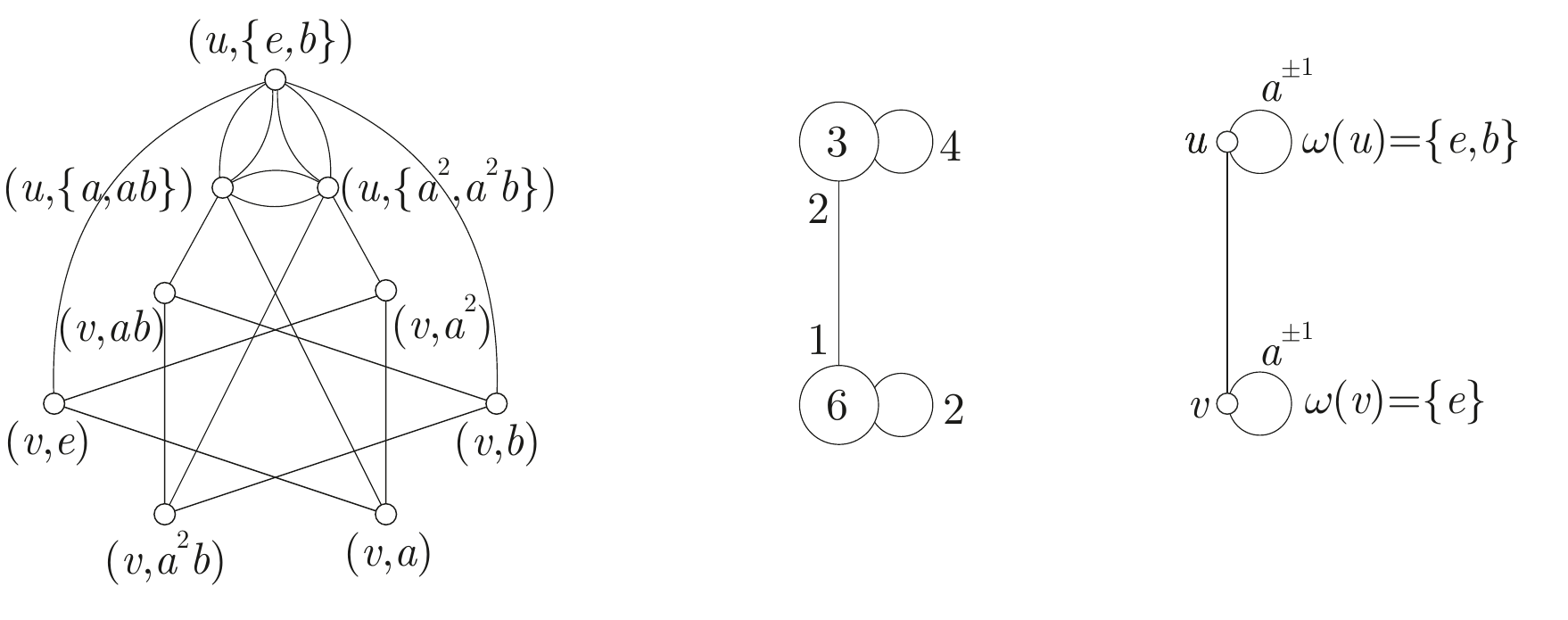}
     \caption{The factored lift graph $\G^{(\alpha,\omega)}$ on the dihedral group $D_3$ of Example \ref{ex:exemple3}, a regular partition of it, and its corresponding combined voltage graph.}
    \label{fig:factored-lift(c)}
\end{figure}
\end{example}

Thus,  $\B$ is the $2\times 2$ matrix shown in \eqref{blockB3}, where $\Sigma(G_u)=e+b$.
\begin{equation}
 \B=\left(
 \begin{array}{ccc}
 (\Sigma(G_u))(a+a^{-1}) & \Sigma(G_u)\\
 e & a+a^{-1}
 \end{array}
 \right).
 \label{blockB3}
\end{equation}
The matrices $\B(\rho_i)$ for $i=0,1,2$ are shown in \eqref{eq:rho01(B)(c)} and \eqref{eq:rho2(B)(c)}:
\begin{equation}
 \B(\rho_0)=\left(
 \begin{array}{cc}
 4 & 2\\
 1 & 2\\
 \end{array}
 \right),\qquad  \B(\rho_1)=\left(
 \begin{array}{cc}
 0 & 0\\
 1 & 2\\
 \end{array}
 \right),
 \label{eq:rho01(B)(c)}
\end{equation}
{\small
\begin{equation}
 \B(\rho_2)=\left(
 \begin{array}{cc|cc}
 \zeta+\zeta^{-1} & \zeta+\zeta^{-1} & 1 & 1\\
 \zeta+\zeta^{-1} & \zeta+\zeta^{-1} & 1 & 1\\
 \hline
 1 & 0 & \zeta+\zeta^{-1} & 0\\
 0 & 1 & 0 & \zeta+\zeta^{-1}
 \end{array}
 \right)=\left(
 \begin{array}{cccc}
 -1 & -1 & 1 & 1\\
 -1 & -1 & 1 & 1\\
 1 & 0 & -1 & 0\\
 0 & 1 & 0 & -1
 \end{array}
 \right).
 \label{eq:rho2(B)(c)}
\end{equation}
}
The corresponding eigenvalues and eigenvectors are listed in Table
\ref{tab:eigenvv-rho01(c)}.
\begin{table}[ht]
\begin{center}
\begin{tabular}{|cc|cc|}
\hline
 $\B(\rho_0)$ &  & $\B(\rho_1)$ & \\
\hline
\hline
  $3+\sqrt{3}$  & $3-\sqrt{3}$  & $2$  &  $0^*$   \\
\hline
$(1{+}\sqrt{3},1)^{\top}$ & $(1{-}\sqrt{3},1)^{\top}$ & $(0,1)^{\top}$ & $(-2,1)^{\top}$ \\
\hline
\end{tabular}
\vskip 10pt
\begin{tabular}{|cccc|}
\hline
   $\B(\rho_2)$   &  &  & \\
\hline
\hline
  $0^+$  & $0^{++}$  & $-1$  & $-3$   \\
\hline
$(1,0,1,0)^{\top}$ & $(0,1,0,1)^{\top}$ & $(0,0,-1,1)^{\top}$ & $(-2,-2,1,1)^{\top}$ \\
\hline
\end{tabular}
\end{center}
\caption{The eigenvalues and eigenvectors of the matrices $\B(\rho_r)$ for $r=0,1,2$ associated with the combined voltage graph of Figure \ref{fig:factored-lift(c)}.}
\label{tab:eigenvv-rho01(c)}
\end{table}

The eigenvector of $\B(\rho_1)$ not satisfying the condition {\bf (C)} corresponds again to the eigenvalue $0$ with an asterisk. But observe that although the eigenvectors of $\B(\rho_2)$ corresponding to the eigenvalues $0^+$ and $0^{++}$ do not satisfy this condition, their sum, that is, the vector $\f^*=(1,1,1,1)^{\top}$ does! This is because $\rho_2(h)\f^*(u)$ and $\rho_2(h)\f^*(v)$ are constant for every $h\in gG_u=gG_w=\{g,gb\}$. Then, by Theorem \ref{theo-sp-d-dim}, the sum eigenvector $\f^*$  produces two eigenvectors in the factored lift that correspond to the eigenvalue $0$. The same reasoning applies to the eigenvector of the eigenvalue $-3$ so that it also satisfies the condition \eqref{eq:Cond-d}. This example was chosen to illustrate the interesting feature that may happen: A sum of two eigenvectors for the eigenvalue $0$ may satisfy the condition {\bf (C)} even if neither of the two vectors does.

Finally, the spectrum of the factored lift $\G^{(\alpha,\omega)}$ is
$$
\spec(\G^{(\alpha,\omega)})=\{\,3+\sqrt{3},\,2,\,3-\sqrt{3},\,0^{[2]},\,-1^{[2]},\,-3^{[2]}\,\}.
$$

To show that the eigenvectors of the factored lift $\Gamma^{(\alpha,\omega)}$ of Figure \ref{fig:factored-lift(c)}
are linearly independent, we apply part $(iii)$ of Theorem \ref{theo-sp-d-dim}. Indeed, such eigenvectors are obtained from the matrix product $\S\T$, where $\S$ is a $12\times 12$ matrix with block form
$$
\S= (\S_0\, |\, \S_1\, |\, \S_2) = \left(
\begin{array}{cc|cc|cccc}
\S_{0,1} &  \O & \S_{1,1} & \O & \S_{2,1} & \O & \S_{2,2} & \O \\
\O & \S_{0,1} &  \O & \S_{1,1} & \O & \S_{2,1} & \O & \S_{2,2}
\end{array}\right),
$$
where
\begin{align*}
\S_{0,1} &=(\rho_0(e),\rho_0(a),\rho_0(a^2),\rho_0(b),\rho_0(ab),\rho_0(a^2b))^{\top}
=(1,1,1,1,1,1)^{\top},\\
\S_{1,1} &=(\rho_1(e),\rho_1(a),\rho_1(a^2),\ldots,\rho_1(a^2b))^{\top}=(1,1,1,-1,-1,-1)^{\top},\\
\S_{2,1} &=(\rho_2(e)_1,\rho_2(a)_1,\rho_2(a^2)_1,\ldots,\rho_2(a^2b)_1)^{\top}
=
\left(
\begin{array}{cccccc}
1 & z & z^2 & 0 & 0 & 0\\
0 & 0 & 0 & 1 & z & z^2
\end{array}
\right)^{\top},
\\
\S_{2,2} &=(\rho_2(e)_2,\rho_2(a)_2,\rho_2(a^2)_2,\ldots,\rho_2(a^2b)_2)^{\top}
=
\left(
\begin{array}{cccccc}
0 & 0 & 0 & 1 & z^{-1} & z^{-2}\\
1 & z^{-1} & z^{-2} & 0 & 0 & 0
\end{array}
\right)^{\top}.
\end{align*}
Note that, as required in the proof, $\S_{2,1}$ and $\S_{2,2}$ are formed, respectively, out of the first and second rows of the matrices $\rho_2(e)$, $\rho_2(a)$, $\rho_2(a^2)$, and so on. The matrix $\T$, in this case, is a $12\times 12$ matrix with block form $\T = \diag(\T_0,\T_1,\T_2)$, where $\T_0 = \U_0$, $\T_1=\U_1$, and $\T_2 = \diag(\U_2,\U_2)$, so that
$$
\T =\left(
\begin{array}{cccc}
\U_0 &  \O & \O & \O \\
\O & \U_1 & \O & \O \\
\O & \O & \U_2 & \O \\
\O & \O & \O & \U_2
\end{array}\right).
$$
Here, one needs to be careful about the indexation of rows and columns to align eigenvectors with the corresponding eigenvalues. In accordance with the proof of Theorem \ref{theo-sp-d-dim}, for each $\rho_i\in\{\rho_0,\rho_1,\rho_2\}$ of dimension $d_i$, the  $d_i\times d_i$ matrix $\U_i$ is formed by a choice of the corresponding eigenvectors of $\rho_i(\B)$. To proceed, we choose to list the eigenvalues in the order $3+\sqrt{3}$, $3-\sqrt{3}$, $2$, $0$, $-1$, $-3$, $0$, $0$, together with a choice of the corresponding eigenvectors as follows:
\begin{equation*}\label{eq:Ui}
\U_0 = \left(
\begin{array}{cc}
1+\sqrt{3}  & 1-\sqrt{3} \\
1  &  1
\end{array}
\right),\quad
\U_1 = \left(
\begin{array}{cc}
0  & -2 \\
1  &  1
\end{array}
\right),\quad 
\U_2 =
\left(
\begin{array}{cccc}
0 & -2  &  1  & 1 \\
0 & -2  &  1  & 0 \\
-1 &  1  &  1  & 1 \\
1  & 1  &  1  &  0
\end{array}\right),
\end{equation*}
where, from left to right, the columns of $\U_0$ correspond to eigenvalues $3+\sqrt{3}$ and $3-\sqrt{3}$, the columns of $\U_1$ to the eigenvalues $2$ and $0$, and finally the columns of $\U_2$ correspond to the eigenvalues $-1, -3,0$ and $0$, respectively. Then, we get
$$
\footnotesize
\S\T=\left(
\begin{array}{cccccccccccc}
1+\sqrt{3} & 1-\sqrt{3} & 0 & -2 & 0 & -2 & 1 & 1 & 0 & -2 & 1 & 0\\
1+\sqrt{3} & 1-\sqrt{3} & 0 & -2 & 0 & -2z & z & z & 0 & -2z^{-1} & z^{-1} & 0\\
1+\sqrt{3} & 1-\sqrt{3} & 0 & -2 & 0 & -2z^2 & z^2 & z^2 & 0 & -2z^{-2} & z^{-2} & 0\\
1+\sqrt{3} & 1-\sqrt{3} & 0 & 2 & 0 & -2 & 1 & 0 & 0 & -2 & 1 & 1\\
1+\sqrt{3} & 1-\sqrt{3} & 0 & 2 & 0 & -2z & z & 0 & 0 & -2z^{-1} & z^{-1} & z^{-1}\\
1+\sqrt{3} & 1-\sqrt{3} & 0 & 2 & 0 & -2z^2 & z^2 & 0 & 0 & -2z^{-2} & z^{-2} & z^{-2}\\
1 & 1 & 1 & 1 & -1 & 1 & 1 & 1 & 1 & 1 & 1 & 0\\
1 & 1 & 1 & 1 & -z & z & z & z & z^{-1} & z^{-1} & z^{-1} & 0\\
1 & 1 & 1 & 1 & -z^2 & z^2 & z^2 & z^2 & z^{-2} & z^{-2} & z^{-2} & 0\\
1 & 1 & -1 & -1 & 1 & 1 & 1 & 0 & -1 & 1 & 1 & 1\\
1 & 1 & -1 & -1 & z & z & z & 0 & -z^{-1} & z^{-1} & z^{-1} & z^{-1}\\
1 & 1 & -1 & -1 & z^2 & z^2 & z^2 & 0 & -z^{-2} & z^{-2} & z^{-2} & z^{-2}
\end{array}
\right).
$$
Then, the eigenvectors of the lift are obtained, first removing the columns $4,8,12$, where the first and fourth rows  (corresponding to the elements of $G_u=\{e,b\}$) are different, and second removing the rows $1,2,3$ (whose entries are equal to the rows $4,5,6$, respectively):
$$
\left(
\begin{array}{ccccccccc}
1+\sqrt{3} & 1-\sqrt{3} & 0 & 0 & -2 & 1 & 0 & -2 & 1 \\
1+\sqrt{3} & 1-\sqrt{3} & 0 & 0 & -2z & z & 0 & -2z^{-1} & z^{-1} \\
1+\sqrt{3} & 1-\sqrt{3} & 0 & 0 & -2z^2 & z^2 & 0 & -2z^{-2} & z^{-2} \\
1 & 1 & 1 & -1 & 1 & 1 & 1 & 1 & 1 \\
1 & 1 & 1 & -z & z & z & z^{-1} & z^{-1} & z^{-1}\\
1 & 1 & 1 & -z^2 & z^2 & z^2 & z^{-2} & z^{-2} & z^{-2} \\
1 & 1 & -1 & 1 & 1 & 1 & -1 & 1 & 1 \\
1 & 1 & -1 & z & z & z & -z^{-1} & z^{-1} & z^{-1} \\
1 & 1 & -1 & z^2 & z^2 & z^2 & -z^{-2} & z^{-2} & z^{-2}
\end{array}
\right).
$$


\section{Concluding remarks}\label{sec:rem}

In Section \ref{sec:factored-lifts}, we presented a method for determining a complete spectrum of a factored lift from the spectrum of a special matrix reflecting the structure of a base graph, together with a combined voltage assignment in a group, with entries in a complex group algebra associated with the voltage group. For a similar derivation of eigenvectors of the lift, we derived a sufficient condition, and in Section \ref{sec:examples}, we illustrated the complexity of the situation with lifting eigenvectors in general.

In a way, this is a bit of a paradox when one considers the previous work \cite{dfps21,dfps23} of the same set of authors. The two papers offer a complete description of permutation lifts of both spectra and eigenspaces, together with details of the underpinning theory in \cite{dfps21}. Permutation lifts represent a generalisation of ordinary lifts in a completely different direction compared with factored lifts (as explained at the end of Section \ref{sec:comb-fact}).

The factored lifts of base graphs equipped with a combined voltage assignment in a given group, considered in this paper, are an equivalent way of studying the quotients of graphs that admit a free action of a given group on arcs. As alluded to earlier, an investigation of quotients of graphs by general subgroups of automorphisms and a formal study of their reconstruction by `general lifts' was initiated by Poto\v{c}nik and Toledo in \cite{pt21}. 


\subsection*{Acknowledgment}
\thanks{The first two authors' research 
was supported by project PID2025-171721NB-I00 funded by MICIU/AEI/10.13039/501100011033 and by ERDF, EU; and by project
2021SGR00434 funded by AGAUR. 
 The third and fourth authors acknowledge support of this research by the APVV Research Grants 22-0005 and 23-0076, and the VEGA Research Grants 1/0069/23 and 1/0011/25.}


\end{document}